\documentclass[10pt,a4paper,english]{endm}
\usepackage[T1]{fontenc}
\usepackage[latin9]{inputenc}
\usepackage{calc}
\usepackage{amsmath}
\usepackage{amssymb}
\usepackage{stmaryrd}

\makeatletter

\usepackage{endmmacro}
\usepackage{graphicx}

\usepackage{tikz,tkz-tab,tkz-graph}
\usepackage{tikz}

\usetikzlibrary{arrows,automata}
\usetikzlibrary{positioning, fadings}
\tikzset{
    state/.style={
           rectangle,
  fill=#1!5!white,
           rounded corners,
           draw=#1, very thick,
           minimum height=2em,
           inner sep=2pt,
           text centered,
           },
coeff/.style={
           circle,
           draw=black, very thick,
           minimum height=2em,
           inner sep=2pt,
           text centered,
           },	
}

\AtBeginDocument{

}

\makeatother

\usepackage{babel}
\begin{document}
\begin{verbatim}\end{verbatim}\vspace{2.5cm}
\begin{frontmatter}
\title{On Adjacency and e-Adjacency in General Hypergraphs: Towards a New e-Adjacency Tensor}
\author[CERN,Unige]{X. Ouvrard\thanksref{SupportCERN}},
\author[CERN]{J.M. Le Goff} and 
\author[Unige]{S. Marchand-Maillet}
\address[CERN]{CERN, CH-1211 Geneva 23}
\address[Unige]{University of Geneva, CUI, 7 route de Drize, Battelle   A, CH-1227 Carouge}
\thanks[SupportCERN]{This document is part of X. Ouvrard PhD work supervised by Pr. S. Marchand-Maillet and J.M. Le Goff and founded by a doctoral position at CERN.}
\thanks[XOMail]{Email: \href{mailto:xavier.ouvrard@cern.ch} {\texttt{\normalshape    xavier.ouvrard@cern.ch}}}
\begin{abstract}
In graphs, the concept of adjacency is clearly defined: it is a pairwise relationship between vertices. Adjacency in hypergraphs has to integrate hyperedge multi-adicity: the concept of adjacency needs to be defined properly by introducing two new concepts: $k$-adjacency -  $k$ vertices are in the same hyperedge - and e-adjacency - vertices of a given hyperedge are e-adjacent. In order to build a new e-adjacency tensor that is interpretable in terms of hypergraph uniformisation, we designed two processes: the first is a hypergraph uniformisation process (HUP) and the second is a polynomial homogeneisation process (PHP). The PHP allows the construction of the e-adjacency tensor while the HUP ensures that the PHP keeps interpretability. This tensor is symmetric and can be fully described by the number of hyperedges; its order is the range of the hypergraph, while extra dimensions allow to capture additional hypergraph structural information including the maximum level of $k$-adjacency of each hyperedge. Some results on spectral analysis are  discussed.
\end{abstract}
\end{frontmatter}

\section{Adjacency in hypergraphs}

A \textbf{hypergraph} $\mathcal{H}=\left(V,E\right)$ is a hyperedge
family $E=\left\{ e_{i}:e_{i}\subseteq V\land i\in\left\llbracket p\right\rrbracket \right\} $\footnote{$\left\llbracket k;n\right\rrbracket $ is $\left\{ i:i\in\mathbb{N}\land k\leqslant i\leqslant n\right\} $
and $\left\llbracket n\right\rrbracket $ is $\left\llbracket 1,n\right\rrbracket $.
$\mathcal{S}_{k}$ is the permutation set on $\left\llbracket k\right\rrbracket $.} over the vertex set $V=\left\{ v_{i}:i\in\left\llbracket n\right\rrbracket \right\} $
\cite{bretto2013hypergraph}. A \textbf{hypergraph with no repeated
hyperedge} is a hypergraph where the hyperedges are distinct pairwise.

We write $k_{\max}=\max\left\{ \left|e\right|:e\in E\right\} $ the
\textbf{range of the hypergraph}. 

Hyperedge multi-adicity calls for additional adjacency concepts.

\begin{definition}

$k$ vertices are said \textbf{$k$-adjacent} if it exists an hyperedge
that contains them. Vertices of a given hyperedge are said \textbf{e-adjacent}.
The $\overline{k}$-adjacency of an hypergraph is the maximal value
of $k$ such that it exists vertices of the hypergraph that are $k$-adjacent.

\end{definition}

Hypermatrices - abusively designated as tensors \cite{qi2017siam}
- are used to store the adjacency multi-adic relationships. In $k$-uniform
hypergraphs, where all hyperedges have the same cardinality $\overline{k}=k$,
$\overline{k}$-adjacency and e-adjacency are equivalent; we use here
the degree normalized $\overline{k}$-adjacency hypermatrix \cite{cooper2012spectra}.

For general hypergraphs with no-repeated hyperedge, a first e-adjacency
hypermatrix is defined in \cite{banerjee2017spectra}. The value and
the number of elements that are required to store this hypermatrix
vary depending on the hyperedge cardinality; due to index repetition,
tensor elements can not be interpreted directly in term of a hypergraph
uniformisation process (HUP). To address this issue, we propose a
new e-adjacency tensor\footnote{Details and proofs can be found in \cite{ouvrard2017cooccurrence}.}. 

\section{A new e-adjacency tensor for general hypergraphs}

We give here only the main steps.\footnote{Exponents into parenthesis refer to the order of the corresponding
tensor; indices into parenthesis refer to a sequence of objects. } 

\subsection{Decomposition in layers}

The family $\left(E_{k}\right)_{1\leqslant k\leqslant k_{\text{\text{max}}}}$
where $E_{k}=\left\{ e\in E:\,\left|e\right|=k\right\} $ constitutes
a partition of $E$. $\mathcal{H}$ is decomposable uniquely into
a $k$-uniform hypergraph direct sum $\mathcal{H}=\bigoplus\limits _{k=1}^{k_{\max}}\mathcal{H}_{k}$
of increasing $k\in\left\llbracket k_{\max}\right\rrbracket $. The
$\mathcal{H}_{k}=\left(V,E_{k}\right)$ - $k\in\left\llbracket k_{\max}\right\rrbracket $
- are called the \textbf{layers} of $\mathcal{H}$. Any of these $\mathcal{H}_{k}$
is representable by a degree-normalised $\overline{k}$-adjacency
hypermatrix $\boldsymbol{A}_{k}=\left(a_{(k)\,i_{1}...i_{k}}\right)$. 

Symmetric cubical hypermatrices are bijectively mapped to homogeneous
polynomials \cite{comon2015polynomial} through the hypermatrix multilinear
matrix multiplication \cite{lim2013tensors}. 

We build a family $P_{\mathcal{H}}=\left(P_{k}\right)$ of homogenous
polynomials that are one-to-one mapped to the layers of the hypergraph.
Considering $\boldsymbol{z}=\left(\boldsymbol{z}_{0}\right)^{\top}$\footnote{We write $\boldsymbol{z}_{0}$ the variable list $z^{1},...,z^{n}$
and $\boldsymbol{z}_{k}$ the variable list $\boldsymbol{z}_{0},y^{1},...,y^{k}$. } - for all $i\in\left\llbracket n\right\rrbracket $: $z^{i}$ represents
$v_{i}\in V$ - and $\left(\boldsymbol{z}\right)_{[k]}=\left(\boldsymbol{z},...,\boldsymbol{z}\right)\in\left(\mathbb{R}^{n}\right)^{k}$,
$\left(\boldsymbol{z}\right){}_{[k]}.\boldsymbol{A_{k}}$ contains
only one element: $P_{k}\left(\boldsymbol{z}_{0}\right)=\sum\limits _{1\leqslant i_{1},...,i_{k}\leqslant n}a_{(k)\,i_{1}...i_{k}}z^{i_{1}}...z^{i_{k}}.$\label{Eq:1}
As $\boldsymbol{A_{k}}$ is symmetric: $P_{k}\left(\boldsymbol{z}_{0}\right)=\sum\limits _{1\leqslant i_{1}\leqslant...\leqslant i_{k}\leqslant n}\alpha_{(k)\,i_{1}...i_{k}}z^{i_{1}}...z^{i_{k}}$
with $\alpha_{(k)\,i_{1}...i_{k}}=k!a_{(k)\,i_{1}...i_{k}}.$

\subsection{Uniformisation and homogeneisation process}

The hypergraph uniformisation process involves two elementary operations
on weighted hypergraphs.

\textbf{\hspace{-1.5em}}%
\noindent\fbox{\begin{minipage}[t]{1\columnwidth - 2\fboxsep - 2\fboxrule}%
\textbf{Operation 1:} Let $\mathcal{H}_{w}=\left(V,E,w\right)$ be
a weighted hypergraph. Let $y\notin V$. 

The \textbf{$y$-vertex-augmented hypergraph} of $\mathcal{H}_{w}$
is the weighted hypergraph $\overline{\mathcal{H}_{\overline{w}}}=\left(\overline{V},\overline{E},\overline{w}\right)$
where $\overline{V}=V\cup\left\{ y\right\} $, $\overline{E}=\left\{ \phi\left(e\right):e\in E\right\} $
- with the map $\phi:\mathcal{P}\left(V\right)\rightarrow\mathcal{P}\left(\overline{V}\right)$
such that: $\forall A\in\mathcal{P}\left(V\right):\,$$\phi(A)=A\cup\left\{ y\right\} $
- and, $\overline{w}$ such that $\forall e\in E$: $\overline{w}\left(\phi(e)\right)=w(e).$%
\end{minipage}}

\textbf{\hspace{-1.5em}}%
\noindent\fbox{\begin{minipage}[t]{1\columnwidth - 2\fboxsep - 2\fboxrule}%
\textbf{Operation 2:} The \textbf{merged hypergraph} $\widehat{\mathcal{H}_{\widehat{w}}}=\left(\widehat{V},\widehat{E},\widehat{w}\right)$
of two weighted hypergraphs $\mathcal{H}_{a}=\left(V_{a},E_{a},w_{a}\right)$
and $\mathcal{H}_{b}=\left(V_{b},E_{b},w_{b}\right)$ is the weighted
hypergraph with vertex set $\widehat{V}=V_{a}\cup V_{b}$, with hyperedge
family $\widehat{E}=E_{a}+E_{b}$ - constituted of all elements of
$E_{a}$ and all elements of $E_{b}$ - such that $\forall e\in E_{a}$,
$\widehat{w}(e)=w_{a}(e)$ and $\forall e\in E_{b}$, $\widehat{w}(e)=w_{b}(e).$%
\end{minipage}}

\textbf{The hypergraph uniformisation process} starts by mapping each
$\mathcal{H}_{k}$ to a weighted hypergraph $\mathcal{H}_{w_{k},k}=\left(V,E_{k},w_{k}\right)$
with: $\forall e\in E_{k}:w_{k}(e)=c_{k}$ with $c_{k}\in\mathbb{R}^{+*}$
and $k\in\left\llbracket k_{\max}\right\rrbracket $. $c_{k}$ are
dilatation coefficients introduced to guarantee that the generalized
hand-shake lemma holds in the e-adjacency tensor. A set of pairwise
distinct vertices $V_{s}=\left\{ y_{k}:k\in\left\llbracket k_{\max}-1\right\rrbracket \right\} $
is generated and such that no vertex of $V_{s}$ is in $V$.

The HUP iterates over a two-phase step: the \textbf{inflation phase}
(IP) and the \textbf{merging phase} (MP). At step $k>1$ the input
is the $\left(k-1\right)$-uniform weigthed hypergraph $\mathcal{K}_{w}$
obtained from the previous iteration; at step 1, $\mathcal{K}_{w}=\mathcal{H}_{w_{1},1}.$
In the IP, $\mathcal{K}_{w}$ is transformed into $\overline{\mathcal{K}_{\overline{w}}}$
the $k+1$-uniform $y_{k}$-vertex-augmented hypergraph of $\mathcal{K}_{w}$.

The MP elaborates the merged hypergraph $\widehat{\mathcal{K}_{\widehat{w}}}$
from $\overline{\mathcal{K}_{\overline{w}}}$ and $\mathcal{H}_{w_{k+1},k+1}.$

At the end of each step $k$ is increased until it reaches $k_{\max}$:
the last $\widehat{\mathcal{H}_{\widehat{w}}}$ obtained is called
the \textbf{$V_{s}$-layered uniform hypergraph} of $\mathcal{H}$.

\begin{proposition}

$\widehat{\mathcal{H}_{\widehat{w}}}$ captures exactly the e-adjacency
of $\mathcal{H}$.

\end{proposition}

\textbf{In the polynomial homogeneisation process}, $R_{\mathcal{H}}=\left(R_{k}\right)_{k\in\left\llbracket k_{\max}\right\rrbracket }$
the family of homogeneous polynomials of degree $k$ is obtained iteratively
from the family $\left(c_{k}P_{k}\right)_{k\in\left\llbracket k_{\max}\right\rrbracket }$:
for all $k\in\left\llbracket k_{\max}\right\rrbracket $, $c_{k}P_{k}$
maps one to one to $\mathcal{H}_{w_{k},k}$. 

We set $R_{1}\left(\boldsymbol{z}_{o}\right)=c_{1}P_{1}\left(\boldsymbol{z}_{o}\right)=c_{1}\sum\limits _{i=1}^{n}a_{(1)\,i}z^{i}.$
We generate $k_{\max}-1$ new pairwise distinct variables $y^{j}$,
$j\in\left\llbracket k_{\max}-1\right\rrbracket $.

At step $k$, we suppose that: $R_{k}\left(\boldsymbol{z}_{k-1}\right)=\sum\limits _{j=1}^{k}c_{j}\sum\limits _{i_{1},...,i_{j}=1}^{n}a_{(j)\,i_{1}...i_{j}}z^{i_{1}}...z^{i_{j}}\prod\limits _{l=j}^{k-1}y^{l},$
with the convention that: $\prod\limits _{l=j}^{k-1}y^{l}=1$ if $j>k-1.$
Then for $y^{k-1}\neq0$:

\begin{eqnarray*}
R_{k+1}\left(\boldsymbol{z}_{k}\right) & = & y^{k\,(k+1)}\left(R_{k}\left(\dfrac{\boldsymbol{z}_{k-1}}{y^{k\,(k)}}\right)+c_{k+1}P_{k+1}\left(\dfrac{\boldsymbol{z}_{o}}{y^{k\,(k+1)}}\right)\right)\\
 & = & R_{k}\left(\boldsymbol{z}_{k-1}\right)y^{k}+c_{k+1}\sum\limits _{i_{1},...,i_{k+1}=1}^{n}a_{(k+1)\,i_{1}\,...\,i_{k+1}}z^{i_{1}}...z^{i_{k+1}}
\end{eqnarray*}

and for $y^{k}=0$: $R_{k+1}\left(\boldsymbol{z}_{k-1},0\right)=c_{k+1}\sum\limits _{i_{1},...,i_{k+1}=1}^{n}a_{(k+1)\,i_{1}\,...\,i_{k+1}}z^{i_{1}}...z^{i_{k+1}}.$

Even if $P_{k+1}\left(z_{0}\right)=0$ the step above is performed:
the degree of $R_{k}$ will increase by 1.

\subsection{Construction of the e-adjacency tensor}

From $R_{\mathcal{H}}=\left(R_{k}\right)$ we build a symmetric tensor.
$R_{k}$ is an homogeneous polynomial with $n+k-1$ variables of order
$k$. With $w_{\left(k\right)}$ for $w_{\left(k\right)}^{1},...,w_{\left(k\right)}^{n}$,
we have: $R_{k}\left(\boldsymbol{w}_{(k)}\right)=\sum\limits _{i_{1},...,i_{k}=1}^{n+k-1}r_{(k)\,i_{1}\,...\,i_{k}}w_{(k)}^{i_{1}}...w_{(k)}^{i_{k}}$
where:
\begin{itemize}
\item for $i\in\left\llbracket n\right\rrbracket $: $w_{(k)}^{i}=z^{i}$
and for $i\in\left\llbracket n+1;n+k-1\right\rrbracket $: $w_{(k)}^{i}=y^{i-n}$
\item for all $\forall j\in\left\llbracket k\right\rrbracket $, for $1\leqslant i_{1}<...<i_{j}\leqslant n$,
for all $l\in\left\llbracket j+1;k\right\rrbracket $\footnote{With the convention $\left\llbracket p,q\right\rrbracket =\emptyset$
if $p>q$}: $i_{l}=n+l-1$ and, for all $\sigma\in\mathcal{S}_{k}$: 
\end{itemize}
\[
r_{(k)\,\sigma\left(i_{1}\right)...\sigma\left(i_{k}\right)}=\dfrac{c_{j}\alpha_{(j)\,i_{1}...i_{j}}}{k!}=\dfrac{j!}{k!}c_{j}a_{(j)\,i_{1}...i_{j}}
\]

\begin{itemize}
\item otherwise $r_{(k)\,i_{1}\,...\,i_{k}}$ is null.
\end{itemize}
Also $R_{k}$ can be linked to a symmetric hypercubic tensor of order
$k$ and dimension $n+k-1$ written $\boldsymbol{R_{k}}$ whose elements
are $r_{(k)\,i_{1}\,...\,i_{k}}$.

The coefficients $c_{k}$, $k\in\left\llbracket k_{\text{max}}\right\rrbracket $
are chosen so that the number of edges calculated by the generalized
handshake lemma is valid. 

We choose: $c_{j}=\dfrac{k_{\max}}{j}$ as:

$\left|E\right|=\dfrac{1}{k_{\max}}\sum\limits _{i_{1},...,i_{k_{\max}}\in\left\llbracket n+k_{\max}-1\right\rrbracket }r_{i_{1}...i_{k_{\max}}}$$=\sum\limits _{j=1}^{k_{\text{\text{max}}}}\dfrac{1}{j}\sum\limits _{i_{1},...,i_{j}\in\left\llbracket n\right\rrbracket }a_{(j)\,i_{1}...i_{j}}.$

Hence, combining above with the fact that $a_{(j)\,i_{1}...i_{j}}=\dfrac{1}{(j-1)!}$
when $\left\{ v_{i_{1}},...,v_{i_{j}}\right\} \in E$ and 0 otherwise:
$r_{i_{1}...i_{k_{\text{\text{max}}}}}=\dfrac{1}{\left(k_{\text{\text{max}}}-1\right)!}$
for nonzero elements of $\boldsymbol{R_{k_{\max}}}$.

\begin{definition}

The hypermatrix $\boldsymbol{R_{k_{\max}}}$ is called the \textbf{layered
e-adjacency tensor} of the hypergraph $\mathcal{H}$. We write it
later $\mathcal{A}_{\mathcal{H}}.$

\end{definition}

\section{Further comments and results}

The HUP adds vertices in the IPs; they give indication on the original
cardinality of the hyperedge they are added to as well as the level
of $k$-adjacency possible in this hyperedge. The resulting tensor
is symmetric and is bijectively associated to the original hypergraph,
containing its overall structure.

We consider in the following propositions a hypergraph $\mathcal{H}=\left(V,E\right)$
with no repeated hyperedge with layered e-adjacency tensor $\mathcal{A}_{\mathcal{H}}=\left(a_{i_{1}...i_{k_{\max}}}\right).$

\begin{proposition}

It holds: $\sum\limits _{\substack{i_{2},...,i_{k_{\max}}=1\\
\delta_{ii_{2}...i_{k_{\max}}=0}
}
}^{n+k_{\max}-1}a_{ii_{2}...i_{k_{\max}}}=d_{i}$

where: $\forall i\in\left\llbracket n\right\rrbracket :\,d_{i}=\deg\left(v_{i}\right)$
and $\forall i\in\left\llbracket k_{\max}-1\right\rrbracket :\,d_{n+i}=\deg\left(y_{i}\right).$

Moreover: $\forall j\in\left\llbracket 2;k_{\max}\right\rrbracket $:
$\left|\left\{ e\,:\,\left|e\right|=j\right\} \right|=d_{n+j}-d_{n+j-1}$ 

and: $\left|\left\{ e\,:\,\left|e\right|=1\right\} \right|=d_{n+1}$

\end{proposition}

Using the definition of eigenvalue of \cite{qi2017siam}, we state:

\begin{theorem}The e-adjacency tensor $\mathcal{A}_{\mathcal{H}}$
has its eigenvalues $\lambda$ such that: 
\begin{equation}
\left|\lambda\right|\leqslant\max\left(\Delta,\Delta^{\star}\right)\label{eq:bound_max_degree_layer}
\end{equation}
 where $\Delta=\underset{1\leqslant i\leqslant n}{\max}\left(d_{i}\right)$
and $\Delta^{\star}=\underset{1\leqslant i\leqslant k_{\max}-1}{\max}\left(d_{n+i}\right).$

\end{theorem}

\begin{proposition}Let $\mathcal{H}$ be a $r$-regular\footnote{A hypergraph is said $r$-regular if all vertices have same degree
$r$.} $r$-uniform hypergraph with no repeated hyperedge. Then this maximum
is reached.

\end{proposition}

\section{Conclusion}

Properly defining the concept of adjacency in a hypergraph is important
to build a proper e-adjacency tensor that preverves the information
on the structure of the hypergraph. The resulting tensor allows to
reconstruct with no ambiguity the original hypergraph. First results
on spectral analysis show that additional vertices inflate the spectral
radius bound. The HUP is a strong basis for further proposals: to
allow repetition of vertices, we introduce hb-graphs, family of multisets
and, propose two other e-adjacency tensors \cite{ouvrard2018adjacency}.

\bibliographystyle{elsarticle-num}
\addcontentsline{toc}{chapter}{\bibname}\bibliography{/home/xo/cernbox/these/000-thesis_corpus/biblio/references}

\begin{thebibliography}{1}
\expandafter\ifx\csname url\endcsname\relax
  \def\url#1{\texttt{#1}}\fi
\expandafter\ifx\csname urlprefix\endcsname\relax\def\urlprefix{URL }\fi
\expandafter\ifx\csname href\endcsname\relax
  \def\href#1#2{#2} \def\path#1{#1}\fi

\bibitem{bretto2013hypergraph}
A.~Bretto, Hypergraph theory, An introduction. Mathematical Engineering. Cham:
  Springer.

\bibitem{qi2017siam}
L.~Qi, Z.~Luo, Tensor analysis: spectral theory and special tensors, Vol. 151,
  SIAM, 2017.

\bibitem{cooper2012spectra}
J.~Cooper, A.~Dutle, Spectra of uniform hypergraphs, Linear Algebra and its
  Applications 436~(9) (2012) 3268--3292.

\bibitem{banerjee2017spectra}
A.~Banerjee, A.~Char, B.~Mondal, Spectra of general hypergraphs, Linear Algebra
  and its Applications 518 (2017) 14--30.

\bibitem{ouvrard2017cooccurrence}
X.~Ouvrard, J.-M. {Le Goff}, S.~Marchand-Maillet, Adjacency and tensor
  representation in general hypergraphs part 1: e-adjacency tensor
  uniformisation using homogeneous polynomials, arXiv preprint
  arXiv:1712.08189.

\bibitem{comon2015polynomial}
P.~Comon, Y.~Qi, K.~Usevich, A polynomial formulation for joint decomposition
  of symmetric tensors of different orders, in: International Conference on
  Latent Variable Analysis and Signal Separation, Springer, 2015, pp. 22--30.

\bibitem{lim2013tensors}
L.-H. Lim, Tensors and hypermatrices, Handbook of Linear Algebra, 2nd Ed., CRC
  Press, Boca Raton, FL (2013) 231--260.

\bibitem{ouvrard2018adjacency}
X.~Ouvrard, J.-M.~L. Goff, S.~Marchand-Maillet, Adjacency and tensor
  representation in general hypergraphs. part 2: Multisets, hb-graphs and
  related e-adjacency tensors, arXiv preprint arXiv:1805.11952.

\end{thebibliography}

\end{document}